\documentclass{article}
\usepackage{graphicx, amsmath, amssymb, amsthm, mathtools, tikz, caption, subcaption, float, natbib}
\title{Domination Density and an Imbalance Regime for Vizing’s Conjecture}
\author{Noah Hosking}
\date{February 2026}
\newtheorem{theorem}{Theorem}
\newtheorem{lemma}{Lemma}
\newtheorem{corollary}{Corollary}
\newtheorem{proposition}{Proposition}
\newtheorem{conjecture}{Conjecture}
\newtheorem{remark}{Remark}
\begin{document}
\maketitle

\begin{abstract}
We develop a domination density framework for studying Vizing’s conjecture $\gamma(G\square H) \ge \gamma(G)\gamma(H)$. Recasting the conjecture in multiplicative density form, we derive a bipartition imbalance sufficient condition for certain graph pairs. For bipartite $G$, we introduce a constructive imbalance amplification argument: if a minimum dominating set of $G$ is sufficiently concentrated on one side of the bipartition relative to $\delta(H)$, then
\[
\gamma(G\square H) + \tau_X |V(H)| \ge \gamma(G)\gamma(H),
\]
where $\tau_X$ depends explicitly on the local domination concentration.  In particular, whenever $G$ is bipartite with $\delta(G) > \delta(H)$, or for when $G$ and $H$ are bipartite with $\delta(G) \ne \delta(H)$ such a concentration must occur. This yields an explicit additive deficit bound within the same domination density regime. We further show that domination-reducing leaf deletions preserve Vizing’s inequality.  Consequently, for bipartite graphs satisfying $\delta(G) > \delta(H)$, the conjecture  reduces to understanding the stability of Vizing’s inequality under domination-neutral  leaf deletions.
\end{abstract}

\section{Introduction}

Let $G=(V(G),E(G))$ be a graph defined by the vertex set $V(G)$ and edge set $E(G)$. Let $\Delta(G)$ be the maximum degree of graph $G$. Let the Cartesian product of two finite simple graphs $G$ and $H$ be denoted $G \square H$. Let $\gamma(G)$ denote the domination number of $G$ (the minimum cardinality of a dominating set $D$ of $G$). Let the domination density $\delta(G)$ be a normalised form of the domination number defined as
\begin{equation}\label{eq: Domination Density}
    \delta(G) := \frac{\gamma(G)}{|V(G)|}.
\end{equation}

For finite, simple, connected bipartite $G$ with bipartition $V(G) = A_G \cup B_G$. Let $D$ be some dominating set of $G$ and let $X \in \{A_G, B_G\}$. Define 
\begin{equation}
    D^\ast(X) :=D\cap X \quad \text{and} \quad \delta^\ast(X): = |D^\ast(X)| / |X|.
\end{equation}

In 1968, Vizing conjectured that for any two finite graphs $G$ and $H$, $\gamma(G\square H) \\ \ge \gamma(G)\gamma(H)$ (\cite{vizing1968some}). Despite extensive study, the conjecture remains open, though various graph classes exist for which it is known to be true. We observe $\gamma(G \square H) \ge \gamma(G)\gamma(H) \iff \delta(G\square H)\ge \delta(G)\delta(H)$ which allows us to obtain a vertex bipartition imbalance sufficient condition for when bipartite $G$ with bipartition $V(G) = A_G \cup B_G$ satisfies $\gamma(G \square H) \ge \gamma(G)\gamma(H)$ based on the ratio $|B_G|/|A_G|$ and maximum degree parameters. 

Our main result is the constructive, if there exists a minimum dominating set of a bipartite graph $G$ such that the dominating vertices are sufficiently concentrated on one side of the bipartition -- relative to $\delta(H)$ -- then an imbalance amplification argument via pendant attachments yields two domination number bounds. In particular, for bipartite $G$ and any $H$ satisfying $\delta(G)>\delta(H)$ or $\delta^\ast(X) > \delta(H)$ we derive
\[
\gamma(G\square H) + \tau_X |V(H)| \ge \gamma(G)\gamma(H)
\]
where $\tau_X:=\min\{ |S| : S \subseteq D^\ast(X) \text{ such that } |S| / |X| > \delta(H) \}$. Additionally, under the same conditions we derive a parameter-free inequality
\[
\gamma(G\square H) + (\delta(H)|X|+1)|V(H)| \ge \gamma(G)\gamma(H).
\]

We note the parameter-free inequality is loose. In particular when compared to the Clark-Suen bound for the domination number of Cartesian products (\cite{clark2000inequality}) we derive
\[
    \gamma(G \square H) \ge\frac{1}{2} \gamma(G)\gamma(H) > \gamma(G)\gamma(H) -  (\delta(H)|X|+1)|V(H)|.
\]
However, since the amplification argument used to obtain
\[
    \gamma(G\square H)+\tau_X|V(H)| \geq \gamma(G)\gamma(H)
\]
proceeds by adjoining pendant vertices, it is natural to ask when these attachments may be reversed while preserving Vizing's inequality. If $\ell$ is a leaf of $G$ such that
\[
    \gamma(G-\ell)=\gamma(G)-1,
\]
then leaf deletion preserves Vizing's inequality, since $|\gamma(G \square H)| - |\gamma((G-\ell)\square H)| \le \gamma(H)$. Thus, the remaining non-trivial case is
\[
    \gamma(G-\ell)=\gamma(G),
\]
where the conjectured lower bound is unchanged even though $\gamma((G-\ell)\square H)$ may decrease. This yields a structural reduction principle for Vizing's conjecture on bipartite graphs. If Vizing's inequality is also preserved under such domination-neutral leaf deletions, then the pendant vertices introduced by the amplification argument may be removed without loss. Consequently, the additive term $\tau_X|V(H)|$ can be eliminated, yielding
\[
    \gamma(G\square H)\geq \gamma(G)\gamma(H)
\]
whenever $G$ admits a minimum dominating set $D$ satisfying $\delta^\ast(X)>\delta(H)$. Hence, within this regime, Vizing's conjecture reduces to the stability of Cartesian-product domination under leaf deletions satisfying $\gamma(G-\ell)=\gamma(G)$.

\section{Bipartite Graphs: Imbalance Regime}\label{sec: Bipartite Imbalance Regime}

In 1968 V.Z. Vizing conjectured the following inequality.
\begin{conjecture}[{\cite{vizing1968some}}]\label{con: Vizing's Conjecture}
    If $G$ and $H$ are any two finite, simple graphs then
    \[
        \gamma(G \square H) \ge \gamma(G)\gamma(H).
    \]
\end{conjecture}

We note an equivalent form the Conjecture.

\begin{proposition}\label{pro: Domination Density Vizing's equivalent}
      If $G$ and $H$ are any two finite, simple graphs then, $$\delta(G \square H) \ge \delta(G)\delta(H) \iff  \gamma(G \square H) \ge \gamma(G)\gamma(H).$$
\end{proposition}
\begin{proof}
Since $|V(G \square H)| = |V(G)||V(H)|$, by the definition of Equation~\ref{eq: Domination Density} $\delta(G \square H)$  $\ge \delta(G)\delta(H)$ is equivalent to
\[
    \frac{\gamma(G \square H)}{|V(G)||V(H)|} \ge \frac{\gamma(G)\gamma(H)}{|V(G)||V(H)|} \iff \gamma(G \square H) \ge \gamma(G)\gamma(H).
\]
\end{proof}

Using this reformulation, since $G \square H$ can be partitioned into $|V(H)|$ induced subgraphs isomorphic to $G$ the vertex set $\{(u,v)\in V(G \square H) : u \in D_G\}$ (symmetrically $\{(u,v)\in V(G \square H) : u \in D_H\}$) is a dominating set of $G \square H$, which implies $\delta(G \square H) \le \delta(G)$ and $\delta(G \square H) \le \delta(H)$. Thus Conjecture~\ref{con: Vizing's Conjecture} may be interpreted as asserting that $\min\{\delta(G), \delta(H)\}$ $\ge \delta(G \square H) \ge \delta(G)\delta(H)$. 

Using this a natural approach becomes studying when the domination density of the Cartesian product is tight to this trivial upper bound. Suppose one can show that $\delta(G \square H) \ge c \cdot \min\{\delta(G),\delta(H)\}$ for some constant $c \in (0,1]$. Then Vizing's inequality holds whenever $c \cdot \min\{\delta(G),\delta(H)\} \ge \delta(G)\delta(H).$ Equivalently, since $\delta(G)\delta(H) = \min\{\delta(G),\delta(H)\}\max\{\delta(G),\delta(H)\},$ it suffices that $c \ge \max\{\delta(G),\delta(H)\}$. In this sense, Vizing's conjecture can be viewed as a statement about the extent to which the domination density of the Cartesian product deviates from its natural upper bound.

Motivated by this approach we use some known domination number and density bounds to derive a sufficient condition for when Conjecture~\ref{con: Vizing's Conjecture} is true for graphs $G$ and $H$.

\begin{lemma}\label{lem: Bipartite Imbalance regime}
Let $G$ be a finite, simple, connected, bipartite graph with  $|V(G)| > 1$ and bipartition $V(G) = A_G \cup B_G$ such that $|A_G| \le |B_G|$. Let $H$ be any finite, simple graph. $G$ and $H$ satisfy Conjecture~\ref{con: Vizing's Conjecture} if
    \[
        \frac{|A_G|+|B_G|}{|A_G|}\ge(\Delta(G)+\Delta(H)+1)\delta(H).
    \]
\end{lemma}
\begin{proof}

Consider the set of vertices $A_G$. Since $G$ is connected, simple, and bipartite, every vertex $b \in B_G$ satisfies $N(b) \subseteq A_G$ and thus has at least one neighbour in $A_G$. $A_G$ is therefore a dominating set of $G$ given every $b$ is dominated by all of its neighbours and all $a \in A_G$ dominates itself. This implies $\gamma(G) \le |A_G|$, considering that $|V(G)|=|A_G|+|B_G|$ by dividing $|A_G|$ by $|A_G|+|B_G|$ we yield
\begin{equation}\label{eq 1 lem: Bipartite Imbalance regime}
    \gamma(G) \le |A_G| \iff \delta(G)  \le \frac{|A_G|}{|A_G|+|B_G|}
\end{equation}

Additionally, as shown by \cite{arnautov1974estimation} and \cite{payan1975nombre} each vertex of $G$ dominates at most $\Delta(G)+1$ vertices implying at least $|V(G)|/(\Delta(G)+1)$ vertices are required to dominate $G$. Noting $\Delta(G\square H)=\Delta(G)+\Delta(H)$ yields
\begin{equation}\label{eq 2 lem: Bipartite Imbalance regime}
    \delta(G \square H)=\frac{\gamma(G\square H)}{|V(G\square H)|} \ge \frac{1}{\Delta(G\square H)+1} = \frac{1}{\Delta(G)+\Delta(H)+1}.
\end{equation}

By sandwiching the bounds of equation~\ref{eq 1 lem: Bipartite Imbalance regime} and \ref{eq 2 lem: Bipartite Imbalance regime} a sufficient condition for $\delta(G\square H)\ge \delta(G)\delta(H)$ is
\[
    \frac{1}{\Delta(G)+\Delta(H)+1} \ge \frac{|A_G|\delta(H)}{|A_G|+|B_G|} \iff \frac{|A_G|+|B_G|}{|A_G|}\ge(\Delta(G)+\Delta(H)+1)\delta(H).
\]
\end{proof}

Using lemma~\ref{lem: Bipartite Imbalance regime}, we build a constructive argument defining an operation which increases graph bipartition unevenness whilst maintaining bounded maximum degree. In doing so we derive new bounding domination number behaviour for $G,H$ when $G$ is bipartite.

\begin{theorem}\label{theo:uneven-or-iterate}
Let $G$ be a finite, simple, connected, bipartite graph with $|V(G)| > 1$ and bipartition $V(G) = A_G \cup B_G$. Let $X \in \{A_G, B_G\}$ and let $D$ be a minimum dominating set of $G$. Let $H$ be any finite, simple graph. If  $\delta^\ast(X)  > \delta(H)$ then
\[
    \gamma(G \square H) + \tau_X |V(H)| \ge \gamma(G)\gamma(H)
\]
where $\tau_X:=\min\{ |S| : S \subseteq D^\ast(X) \text{ such that } |S| / |X| > \delta(H) \}$.
\end{theorem}

\begin{proof}
Fix $X,Y\in\{A_G,B_G\}$ such that $X\ne Y$ and $\delta^\ast(X)>\delta(H)$. Let $S\subseteq D^\ast(X)$ such that $\tau_X:=\min\{ |S| : S \subseteq D^\ast(X) \text{ such that } |S| / |X| > \delta(H) \}$. Construct a sequence of graphs $\{G^{(n)}\}$ with $G^{(0)}=G$ such that for each iteration $n \to n+1$ exactly one pendent leaf is attached to every vertex $s\in S$. Per iteration $X$ remains constant, $|Y|$ increases by $\tau_X$ and each vertex in $X$ gains at most one new neighbour. For $G^{(t)}$ with $t$ sufficiently large we have $|X|<|Y|$; without loss of generality we relabel the bipartition so that $X=A_{G^{(t)}}$ and $Y=B_{G^{(t)}}$.

By Lemma~\ref{lem: Bipartite Imbalance regime}, a sufficient condition for Conjecture~\ref{con: Vizing's Conjecture} to hold for graphs $G^{(t)},H$ is
\begin{equation}\label{eq: 1 theo:uneven-or-iterate}
    \frac{|A_{G^{(t)}}|+|B_{G^{(t)}}|}{|A_{G^{(t)}}|} \ge (\Delta(G^{(t)})+\Delta(H)+1)\delta(H).
\end{equation}

Note each iteration $G^{(n)} \to G^{(n+1)}$ increases the left-hand side of Equation~\ref{eq: 1 theo:uneven-or-iterate} by exactly $\tau_X/|X|$ since $|B_{G^{(n+1)}}|=|B_{G^{(n)}}|+\tau_X$ and $|A_{G^{(n+1)}}|=|X|$. Moreover, since each vertex gains at most one new neighbour per iteration, we have $\Delta(G^{(n+1)})\le \Delta(G^{(n)})+1$, so the right-hand side of Equation~\ref{eq: 1 theo:uneven-or-iterate} increases by at most $\delta(H)$. By assumption $\tau_X/|X|>\delta(H)$ meaning the left-hand side grows strictly faster than the right-hand side. Hence there exists a constant $C$ such that $\gamma(G^{(t)}\square H)\ge \gamma(G^{(t)})\gamma(H)$ for all $t\ge C$.

Attaching a leaf to a vertex $w\in D$ does not increase the domination number: the original dominating set $D$ still dominates all original vertices, and each new leaf is dominated by its neighbour in $D$, hence $\gamma(G^{(t)})=\gamma(G)$. 

Each leaf attachment contributes $|V(H)|$ new vertices to the Cartesian product $G^{(t)}\square H$. Let $Z$ be a minimum dominating set of $G\square H$, viewed as a subset of $V(G^{(t)}\square H)$ on the original fibres. By including all vertices $(s,u)$ with $s\in S$ and $u\in V(H)$, every vertex $(\ell,u)$ in each new leaf fibre $\ell\times H$ is dominated (since $(\ell,u)$ is adjacent to $(s,u)$). 

Therefore $\gamma(G^{(t)}\square H)\le |Z|+\tau_X|V(H)|=\gamma(G\square H)+\tau_X|V(H)|$. Combining inequalities yields $\gamma(G\square H)+\tau_X|V(H)| \ge \gamma(G^{(t)}\square H) \ge \gamma(G)\gamma(H)$ which yields the claim. The minimal choice of $S$ yields the sharpest bound.
\end{proof}

By noting the minimality of $\tau_X$ we derive the following parameter free conditional inequality from Theorem~\ref{theo:uneven-or-iterate}:

\begin{corollary}\label{cor: new bound without tau}
Let $G$ be a finite, simple, connected, bipartite graph with $|V(G)| > 1$ and bipartition $V(G) = A_G \cup B_G$. Let $X \in \{A_G, B_G\}$ and let $D$ be a minimum dominating set of $G$. Let $H$ be any finite, simple graph. If  $\delta^\ast(X)  > \delta(H)$ then
\[
    \gamma(G \square H) + (\delta(H)|X|+1)|V(H)| \ge \gamma(G)\gamma(H).
\]
\end{corollary}
\begin{proof}
Fix $X\in\{A_G,B_G\}$ such that $\delta^\ast(X)>\delta(H)$, and let $\tau_X=\min\big\{|S|: S\subseteq D^\ast(X)\text{ and } |S|/|X|>\delta(H)\big\}$. Since $G$ is connected of order greater than one we have $|X|>0$. By the minimality of $\tau_X$, the integer $\tau_X-1$ does not satisfy the defining inequality, hence
\[
\frac{\tau_X-1}{|X|}\le \delta(H).
\]
Therefore $\tau_X\le \delta(H)|X|+1$. Applying Theorem~\ref{theo:uneven-or-iterate} gives $\gamma(G\square H)+(\delta(H)|X|+1)|V(H)| \ge \gamma(G\square H)+\tau_X|V(H)| \ge \gamma(G)\gamma(H)$.
\end{proof}

Given the condition of Theorem~\ref{theo:uneven-or-iterate} requires just one dominating set of $G$ to satisfy $\delta^\ast(X)  > \delta(H)$ we have two new bounds which very often will be applicable, though not universally so. Thus it becomes relevant to explore conditions for when the inequality $\delta^\ast(X)  > \delta(H)$ is guaranteed for a graph pair. To begin, we consider the implications of assuming $\delta(G) > \delta(H)$.

\begin{corollary}\label{cor: Larger domination density corollary}
Let $G$ be a finite, simple, connected, bipartite graph with $|V(G)| > 1$ and bipartition $V(G) = A_G \cup B_G$. Let $H$ be any finite, simple graph. If $\delta(G)>\delta(H)$ then there exists a minimum dominating set $D \subseteq V(G)$ such that  $\delta^\ast(X)  > \delta(H)$ for some $X \in \{A_G, B_G\}$.
\end{corollary}
\begin{proof}
Let $D$ be a minimum dominating set of $G$, and define $D^\ast(X)=D\cap X$ and $\delta^\ast(X)=|D^\ast(X)|/|X|$ for some $X\in\{A_G,B_G\}$. Assume for contradiction that the hypothesis of Theorem~\ref{theo:uneven-or-iterate} fails for both sides of the bipartition, i.e.
\[
\delta^\ast(A_G)\le \delta(H)
\quad\text{and}\quad
\delta^\ast(B_G)\le \delta(H).
\]
Since $\delta(G)$ is a weighted average of $\delta^\ast(A_G)$ and $\delta^\ast(B_G)$ such that
\[
\delta(G)
=\frac{|D^\ast(A_G)|+|D^\ast(B_G)|}{|A_G|+|B_G|},
\]
it follows that $\delta(G)\le \delta(H)$ which contradicts the assumption $\delta(G)>\delta(H)$. Hence $\delta^\ast(X)>\delta(H)$ for one $X\in\{A_G,B_G\}$.
\end{proof}

If we assume both $G$ and $B$ are bipartite we derive a stronger result but note specific conditions under which Theorem~\ref{theo:uneven-or-iterate} does not apply.

\begin{corollary}\label{cor: Bipartition Local Domination 2}
Let $G$ and $H$ be finite, simple, connected, bipartite graphs such that $|V(G)| > 1$, $|V(H)|>1$, with bipartitions $V(G) = A_G \cup B_G$, and $V(H) = A_H \cup B_H$. If $\delta(G) \ne \delta(H)$ then there exists $\delta^\ast(X)  > \delta(H)$, $X \in \{A_G, B_G\}$ or $\delta^\ast(Y)  > \delta(G)$, $Y \in \{A_H, B_H\}$.
\end{corollary}
\begin{proof}
Assume $G$ and $H$ be finite, simple, connected, bipartite graphs such that $|V(G)| > 1$, $|V(H)|>1$, with bipartitions $V(G) = A_G \cup B_G$, and $V(H) = A_H \cup B_H$. Assuming $\delta(G) \ne \delta(H)$ implies $\delta(G) > \delta(H)$ or $\delta(H) > \delta(G)$, in either case the conditions for corollary~\ref{cor: Larger domination density corollary} are satisfied.
\end{proof}

We also consider the simple restricted case of when a bipartite graph $G$ satisfies $\gamma(G)=|A_G|$. These graphs are considered to be extreme cases for when there exists a subset of dominating vertices in $G$ such that $\delta^\ast(X)  > \delta(H)$ for all $H$ with $|V(H)| >1$:

\begin{corollary}\label{cor: extreme-proposition}
Let $G$ be a finite, simple, connected, bipartite graph with $|V(G)| > 1$ and bipartition $V(G) = A_G \cup B_G$. Let $H$ be any finite, simple, connected graph with $|V(H)| > 1$. If $\gamma(G) = \min \{|A_G|,|B_G| \}$ then there exists a minimum dominating set $D \subseteq V(G)$ such that  $\delta^\ast(X)  > \delta(H)$.
\end{corollary}
\begin{proof}
Since $G$ is connected and $\gamma(G) = |A_G|$ We note $A_G$ is a minimum dominating set of $G$. If we choose $X=A_G$ we obtain $\delta^\ast(A_G)=|A_G| / |A_G| = 1$. Since $H$ is connected and $|V(H)|>1$, it has at least one edge, hence $\gamma(H)\le |V(H)|-1$ and therefore $\delta(H)=\gamma(H)/|V(H)|<1 = \delta^\ast(X)$ for some $D \subseteq V(G)$ yielding the claim.
\end{proof}

Given within this regime the local density condition $\delta^\ast(X)  > \delta(H)$ holds uniformly for all graphs $H$. We note that the regime may have utility not only as an endpoint case, but also as a target toward which other graphs G might be transformed. In particular, developing operations that shift a given graphs minimum dominating sets closer to this regime while preserving the relevant domination properties may provide a route to extending the present framework. We leave this direction for future work.

\section{Comparison to Known Bounds}

Here we briefly compare the additive bound of Corollary~\ref{cor: new bound without tau} with classical multiplicative estimates for $\gamma(G \square H)$. A well-known result of Clark and Suen in \cite{clark2000inequality} states:

\begin{theorem}[Clark-Suen~\cite{clark2000inequality}]
For any finite simple graphs $G$ and $H$,
\[
\gamma(G \square H) \ge \frac{1}{2}\,\gamma(G)\gamma(H).
\]
\end{theorem}

For connected bipartite $G$ with $|V(G)| >1$ we know $\gamma(G) \le \min\{|A_G|,|B_G|\} \le |X|$. Consequently the additive term of Corollary~\ref{cor: new bound without tau}  satisfies $(\delta(H)|X|+1)|V(H)|=\gamma(H)|X|+|V(H)| > \gamma(G)\gamma(H)$. Thus it must hold that
\[
    \gamma(G \square H) \ge\frac{1}{2} \gamma(G)\gamma(H) > \gamma(G)\gamma(H) -  (\delta(H)|X|+1)|V(H)|,
\]
implying the Clark-Suen bound is tighter than Corollary~\ref{cor: new bound without tau}. Accordingly, the additive inequality of Corollary~\ref{cor: new bound without tau}  is not intended to act as a refinement of the Clark-Suen bound. Rather its significance lies in isolating a structural regime: under the hypothesis $\delta^\ast(X)>\delta(H)$, the parameter $\tau_X$ reflects the maximal bounded increase in domination arising from the leaf-attachment construction. Given this, it becomes relevant to analyse the inverse operation to the leaf-attachment step of Theorem~\ref{theo:uneven-or-iterate}: leaf deletions with $\gamma(G-\ell)=\gamma(G)$ with the aim being to minimise $\tau_X$.

\section{Leaf Deletions and Stability of Vizing's Inequality}

The pendant vertex additions in the proof of Theorem~\ref{theo:uneven-or-iterate} preserve $\gamma(G)$ but introduces additional graph fibres in $G^{(t)} \square H$. Since the term $\tau_X |V(H)|$ arises solely from this leaf-attachment step, it's natural to examine the inverse operation and ask whether the removal of pendant vertices from a graph satisfies certain domination properties.

Let $\ell$ be a leaf of $G$ with unique neighbour $s$, and let
\[
    G' = G-\ell.
\]
For any dominating set $D$ of $G\square H$, replace each vertex $(\ell,v)\in D$ by $(s,v)$. The resulting set dominates $G'\square H$ and has cardinality at most $|D|$. Hence
\[
    \gamma(G'\square H)\leq \gamma(G\square H).
\]
Conversely, if $D'$ is a minimum dominating set of $G'\square H$ and $D_H$ is a minimum dominating set of $H$, then $D' \cup \bigl(\{\ell\}\times D_H\bigr)$ dominates $G\square H$. Therefore $\gamma(G\square H) \leq \gamma(G'\square H)+\gamma(H).$ Defining
\[
    \varepsilon
    :=\gamma(G\square H)-\gamma(G'\square H),
\]
we consequently have
\begin{equation}\label{eq: leaf product drop}
    0\leq \varepsilon\leq \gamma(H).
\end{equation}

There are two possible effects of deleting a leaf on the domination number of the first factor. First suppose
\[
    \gamma(G')=\gamma(G)-1.
\]
If $G$ and $H$ satisfy Vizing's inequality, then by Equation~\ref{eq: leaf product drop},
\[
\begin{aligned}
    \gamma(G'\square H) =\gamma(G\square H)-\varepsilon \geq \gamma(G)\gamma(H)-\gamma(H) =\gamma(G')\gamma(H).
\end{aligned}
\]
Thus Vizing's inequality is preserved under every domination-reducing leaf deletion.

The remaining case is therefore $\gamma(G')=\gamma(G)$, which we refer to as a domination-neutral leaf deletion. In this case the conjectured lower bound is unchanged, while the domination number of the Cartesian product may decrease by $\varepsilon$. To quantify the available margin, define the Vizing surplus
\[
    \sigma(G,H)
    :=\gamma(G\square H)-\gamma(G)\gamma(H).
\]
Since $\gamma(G')=\gamma(G)$,
\[
\begin{aligned}
    \gamma(G'\square H)-\gamma(G')\gamma(H)
    &=\gamma(G\square H)-\varepsilon
      -\gamma(G)\gamma(H)\\
    &=\sigma(G,H)-\varepsilon.
\end{aligned}
\]
Consequently, whenever $G$ and $H$ satisfy Vizing's inequality, the
inequality is preserved under the domination-neutral deletion
$G\to G-\ell$ precisely when
\begin{equation}\label{eq: neutral deletion surplus condition}
    \varepsilon\leq \sigma(G,H).
\end{equation}

This idea motivates a broader principle.

\begin{remark}[Leaf-deletion reduction principle]
\label{rem: Conditional strengthening via pendant deletion}
Suppose Vizing's inequality is preserved under leaf deletion; that is, whenever
$\ell$ is a leaf of a graph $G$,
\[
    \gamma(G\square H)\geq\gamma(G)\gamma(H)
    \quad\Longrightarrow\quad
    \gamma((G-\ell)\square H)\geq\gamma(G-\ell)\gamma(H).
\]
Then the pendant attachments introduced in
Theorem~\ref{theo:uneven-or-iterate} may be successively removed while
preserving Vizing's inequality. Consequently,
\[
    \gamma(G\square H)\geq\gamma(G)\gamma(H)
\]
whenever $G$ admits a minimum dominating set $D$ and a bipartition class
$X\in\{A_G,B_G\}$ satisfying
\[
    \delta^\ast(X)>\delta(H).
\]
In particular, by Corollary~\ref{cor: Bipartition Local Domination 2}, this
would establish Vizing's conjecture for bipartite graph pairs satisfying
$\delta(G)\neq\delta(H)$.
\end{remark}

\section{Conclusion}

We have developed a domination-density approach to Vizing's conjecture for bipartite graphs. By exploiting imbalance between the bipartition classes, we obtained sufficient conditions for Vizing's inequality and introduced a pendant-attachment construction yielding the additive bound
\[
    \gamma(G\square H)+\tau_X|V(H)|\geq\gamma(G)\gamma(H)
\]
whenever a minimum dominating set of $G$ satisfies $\delta^\ast(X)>\delta(H)$.

The inverse leaf-deletion problem provides a natural route to strengthening this result. Domination-reducing leaf deletions preserve Vizing's inequality, leaving only domination-neutral deletions, for which $\gamma(G-\ell)=\gamma(G)$, as the essential obstruction. Consequently, stability of Vizing's inequality under these deletions would eliminate the additive term and convert the imbalance regime into a direct sufficient condition for the conjecture. Understanding this remaining stability problem therefore appears to be a natural direction for further investigation.

\bibliographystyle{plainnat}
\bibliography{bib}

@article{vizing1968some,
  title={Some unsolved problems in graph theory},
  author={Vizing, Vadim G},
  journal={Russian Mathematical Surveys},
  volume={23},
  number={6},
  pages={125--141},
  year={1968}
}

@article{clark2000inequality,
  title={Inequality related to Vizing's conjecture},
  author={Clark, W Edwin and Suen, Stephen},
  journal={the electronic journal of combinatorics},
  volume={7},
  pages={N4--N4},
  year={2000}
}

@article{arnautov1974estimation,
  title={Estimation of the exterior stability number of a graph by means of the minimal degree of the vertices},
  author={Arnautov, Vladimir I},
  journal={Prikl. Mat. i Programmirovanie},
  volume={11},
  number={3-8},
  pages={126},
  year={1974}
}

@article{payan1975nombre,
  title={SUR LE NOMBRE D'ABSORPTION D'UN GRAPHE SIMPLE.},
  author={Payan, Charles},
  year={1975}
}

\end{document}